\documentclass[12pt]{article}
\usepackage{amsmath}
\usepackage{amssymb}
\usepackage{vatola}

\def\q{\quad}
\def\qq{\qquad}
\def\mod{\pmod}
\def\t{\text}
\def\f{\frac}
\def\e{\equiv}

\def\a{\alpha}
\def\qtq#1{\q\t{#1}\q}
\def\sls#1#2{(\f{#1}{#2})}
 
\def\Ls#1#2{\Big(\f{#1}{#2}\Big)}
\let \pro=\proclaim
\let \endpro=\endproclaim

\begin{document}

\par\q\par\q
\centerline {\bf On the number of representations of $n$ as a}
\centerline {\bf \qq \qq linear combination of four triangular
numbers II}

$$\q$$
\centerline{Min Wang$^1$ and Zhi-Hong Sun$^2$}
\par\q\newline
\centerline{$\ ^1$School of Mathematical Sciences, Soochow
University,}
 \centerline{Suzhou, Jiangsu 215006,
 P.R. China}
\centerline{Email: 20144207002@stu.suda.edu.cn} \centerline{$\
^2$School of Mathematical Sciences, Huaiyin Normal University,}
\centerline{Huaian, Jiangsu 223001, P.R. China} \centerline{Email:
zhihongsun@yahoo.com} \centerline{Homepage:
http://www.hytc.edu.cn/xsjl/szh}
\par\q\newline

 \abstract{Let $\Bbb Z$ and $\Bbb N$ be the set of integers
 and the set of positive integers, respectively. For
 $a,b,c,d,n\in\Bbb N$ let $N(a,b,c,d;n)$ be the number of
 representations of $n$ by $ax^2+by^2+cz^2+dw^2$, and
 let $t(a,b,c,d;n)$ be the number of
 representations of $n$ by $ax(x-1)/2+by(y-1)/2+cz(z-1)/2
 +dw(w-1)/2$ $(x,y,z,w\in\Bbb Z$). In this paper
 we reveal the connections
between $t(a,b,c,d;n)$ and $N(a,b,c,d;n)$. Suppose $a,n\in\Bbb N$
and $2\nmid a$. We show that
$$t(a,b,c,d;n)=\frac
23N(a,b,c,d;8n+a+b+c+d)-2N(a,b,c,d;2n+(a+b+c+d)/4)$$ for $(a,b,c,d)=
(a,a,2a,8m),\ (a,3a,8k+2,8m+6),\ (a,3a,8m+4,8m+4)\ (n\equiv
m+\frac{a-1}2 \pmod 2)$ and $(a,3a,16k+4,16m+4)\ (n\equiv
\frac{a-1}2\pmod 2)$. We also obtain explicit formulas for
$t(a,b,c,d;n)$ in the cases $(a,b,c,d)=(1,1,2,8),\
(1,1,2,16),(1,2,3,6),\ (1,3,4,12),\ (1,1,$ $3,4),\ (1,1,5,5),\
(1,5,5,5),\ (1,3,3,12),\ (1,1,1,12),\ (1,1,3,12)$ and $(1,3,3,4)$.
 \par\q
 }
 \endabstract
\let\thefootnote\relax \footnotetext {The second author
is supported by the National Natural Science Foundation of China
(grant No. 11371163).}

\section*{1. Introduction}
\par\q  Let $\Bbb Z$ and $\Bbb N$ be the set of integers
 and the set of positive integers, respectively.
For  $n \in \Bbb N$ let
$$\sigma(n)=\sum_{d \mid n,d\in\Bbb N}d.$$ For convenience
 we define $\sigma(n)=0$ for $n\notin \Bbb N$.
 Let $\Bbb Z^4=\Bbb Z\times \Bbb Z\times \Bbb
Z\times \Bbb Z$. For $a,b,c,d\in\Bbb N$ and $n\in\Bbb N \cup \{0\}$
set
$$N(a,b,c,d;n)=\big|\{(x,y,z,w)\in \Bbb Z^4\ |\ n=ax^2+by^2+cz^2+dw^2
\}\big|$$ and $$t(a,b,c,d;n)=\Big|\Big\{(x,y,z,w)\in \Bbb Z^4\ |\ n\
=a\f{x(x-1)}2+ b\f{y(y-1)}2+c\f{z(z-1)}2+d\f{w(w-1)}2\Big\}\Big|.$$
The numbers $\f{x(x-1)}2\ (x\in\Bbb Z)$ are called triangular
numbers.
\par In 1828 Jacobi showed that
$$N(1,1,1,1;n)=8\sum_{d\mid n,4\nmid d}d.\tag 1.1$$
In 1847 Eisenstein (see [D]) gave formulas for $N(1,1,1,3;n)$ and
$N(1,1,1,5;n)$. From 1859 to 1866 Liouville made about 90
conjectures on $N(a,b,c,d;n)$ in a series of papers. Most
conjectures of Liouville have been proved. See [A1, A2, AALW1-AALW5,
AAW], Cooper's survey paper [C], Dickson's historical comments [D]
and Williams' book [W2].
\par
 Let
$$t'(a,b,c,d;n)=\Big|\Big\{(x,y,z,w)\in \Bbb N^4\ |\ n=a\f{x(x-1)}2+
b\f{y(y-1)}2+c\f{z(z-1)}2+d\f{w(w-1)}2\Big\}\Big|.$$ As $\f12
x(x-1)=\f12(-x+1)(-x)$ we have
$$t(a,b,c,d;n)=16t'(a,b,c,d;n).$$
In [L] Legendre stated that
$$t'(1,1,1,1;n)=\sigma(2n+1).\tag 1.2$$
In 2003, Williams [W1] showed that
$$t'(1,1,2,2;n)=\f 14\sum_{d\mid 4n+3}\big(d-(-1)^{\f{d-1}2}\big).$$
 For $a,b,c,d\in\Bbb N$ with $5\le a+b+c+d\le 8$ let
$$C(a,b,c,d)=16+4i_1(i_1-1)i_2+8i_1i_3,$$
where $i_j$ is the number of elements in $\{a,b,c,d\}$ which are
equal to $j$. When $5\le a+b+c+d\le 7$, in 2005 Adiga, Cooper and
Han [ACH] showed that
$$C(a,b,c,d)t'(a,b,c,d;n)=N(a,b,c,d;8n+a+b+c+d).\tag 1.3$$ When
$a+b+c+d=8$, in 2008 Baruah, Cooper and Han [BCH] proved that
$$C(a,b,c,d)t'(a,b,c,d;n)=N(a,b,c,d;8n+8)-N(a,b,c,d;2n+2).\tag 1.4$$
 In 2009,
Cooper [C] determined $t'(a,b,c,d;n)$ for $(a,b,c,d)=(1,1,1,3),\
(1,3,3,3),$ $(1,2,2,3),\ (1,3,6,6),\ (1,3,4,4),\ (1,1,2,6)$ and
$(1,3,12,12)$.
\par In a previous paper [WS], the authors obtained
explicit
 formulas for $t(a,b,c,d;n)$ in the cases
 $(a,b,c,d)=(1,2,2,4),\ (1,2,4,4),\ (1,1,4,4),\ (1,4,4,4)$,
 $(1,3,3,9)$,
 $(1,1,9,9),\ (1,9,9,$ $9)$, $(1,1,1,9)$, $(1,3,9,9)$ and $(1,1,3,9).$

\par Ramanujan's theta functions $\varphi(q)$ and $\psi(q)$ (see [Be])
are defined by
$$\varphi(q)=\sum_{n=-\infty}^{\infty}q^{n^2}=1+2\sum_{n=1}^{\infty}
q^{n^2}\qtq{and} \psi(q)=\sum_{n=0}^{\infty}q^{n(n+1)/2}\
(|q|<1).\tag 1.5$$ It is evident that for $|q|<1$,
$$\sum_{n=0}^{\infty}N(a,b,c,d;n)q^{n}=\varphi(q^a)
\varphi(q^b)\varphi(q^c)\varphi(q^d),$$
$$\sum_{n=0}^{\infty}t'(a,b,c,d;n)q^{n}=\psi(q^a)\psi(q^b)
\psi(q^c)\psi(q^d).$$ From [BCH, Lemma 4.1] we know that for
$|q|<1$,
 $$\align &\varphi(q)=\varphi(q^4)+2q\psi(q^8),\tag 1.6
 \\&\psi(q)\psi(q^3)=\varphi(q^6)\psi(q^4)+q\psi(q^{12})
 \varphi(q^2),\tag
 1.7
 \\&\psi(q)^2=\varphi(q)\psi(q^2).\tag
 1.8\endalign$$
 By (1.6), for $k\in\Bbb N$,
$$\varphi(q^k)=\varphi(q^{4k})+2q^k\psi(q^{8k})
=\varphi(q^{16k})+2q^{4k}\psi(q^{32k})+2q^k\psi(q^{8k}).\tag 1.9$$
\par In this paper, by using Ramanujan's theta functions we reveal
some connections between $t(a,b,c,d;n)$ and $N(a,b,c,d;n)$. Suppose
$a,n\in\Bbb N$ and $2\nmid a$. We show that
$$t(a,b,c,d;n)=\f
23N(a,b,c,d;8n+a+b+c+d)-2N(a,b,c,d;2n+(a+b+c+d)/4)$$ for $(a,b,c,d)=
(a,a,2a,8m),\ (a,3a,8k+2,8m+6),\ (a,3a,8m+4,8m+4)\ (n\e m+\f{a-1}2
\mod 2)$ and $(a,3a,16k+4,16m+4)\ (n\e \f{a-1}2\mod 2)$. Using the
formulas for $N(a,b,c,d;n)$ in [AALW1-AALW5] and [AAW] we also
obtain explicit formulas for $t(a,b,c,d;n)$ in the cases
$(a,b,c,d)=(1,1,2,8),\ (1,1,2,16),(1,2,3,6),\ (1,3,4,$ $12)$,\
$(1,1,3,4)$, $\ (1,1,5,5),\ (1,5,5,5),\ (1,3,3,12),\ (1,1,1,12),\
(1,1,3,12)$ and $(1,3,3,4)$.
\par Throughout this paper  $\sls am$ is the
 Legendre-Jacobi-Kronecker symbol. For $n\in\Bbb N$,  $a(n)$ is given by
$$q\prod_{n=1}^{\infty}(1-q^{2n})(1-q^{4n})(1-q^{6n})(1-q^{12n})=\sum_{n=1}^{\infty}
 a(n)q^n\ (|q|<1).$$

\section*{2. Four relations between $t(a,b,c,d;n)$ and
$N(a,b,c,d;n)$}

\pro{Theorem 2.1} Let $m,n\in\Bbb N$ and $a\in\{1,3,5,\ldots\}$.
Then
$$t(a,a,2a,8m;n)=\f 23N(a,a,2a,8m;8n+8m+4a)-2N(a,a,2a,8m;2n+2m+a).$$
\endpro
Proof. Suppose $|q|<1$ and $a=2s+1$. By (1.9),
 $$\aligned&\sum_{n=0}^{\infty}N(2s+1,2s+1,4s+2,8m;n)q^{n}
 \\&=\varphi(q^{2s+1})^2\varphi(q^{4s+2})\varphi(q^{8m})
 \\&=\big(\varphi(q^{32s+16})+2q^{8s+4}\psi(q^{64s+32})+2q^{2s+1}\psi(q^{16s+8})\big)^2
\\&\q\times\big(\varphi(q^{16s+8})+2q^{4s+2}\psi(q^{32s+16})\big)\cdot\big(\varphi(q^{32m})+2q^{8m}\psi(q^{64m})\big)
 \\&=
\big(\varphi(q^{32s+16})^2+4q^{16s+8}\psi(q^{64s+32})^2
+4q^{4s+2}\psi(q^{16s+8})^2
\\&\q+4q^{8s+4}\varphi(q^{32s+16})\psi(q^{64s+32})
+4q^{2s+1}\varphi(q^{32s+16}))\psi(q^{16s+8})
\\&\q+8q^{10s+5}\psi(q^{64s+32})\psi(q^{16s+8})\big)
\\&\q\times
\big(\varphi(q^{16s+8})\varphi(q^{32m})+2q^{8m}\varphi(q^{16s+8})
\psi(q^{64m})+2q^{4s+2}\psi(q^{32s+16})\varphi(q^{32m})
\\&\q+4q^{8m+4s+2}\psi(q^{32s+16})\psi(q^{64m})\big).
\endaligned\tag 2.1$$
Note that $\varphi(q^{8k_1})^{m_1}\psi(q^{8k_2})^{m_2}
=\sum_{n=0}^{\infty}b_nq^{8n}$ for  any nonnegative integers
$k_1,k_2,m_1$ and $m_2$. From (2.1) we deduce that
$$\align&\sum_{n=0}^{\infty}N(2s+1,2s+1,4s+2,8m;8n+4)q^{8n+4}
\\&=4q^{4s+2}\psi(q^{16s+8})^2\cdot2q^{4s+2}\psi(q^{32s+16})\varphi(q^{32m})
\\&\q+4q^{4s+2}\psi(q^{16s+8})^2\cdot4q^{8m+4s+2}\psi(q^{32s+16})\psi(q^{64m})
\\&\q+4q^{8s+4}\varphi(q^{32s+16})\psi(q^{64s+32})\cdot\varphi(q^{16s+8})\varphi(q^{32m})
\\&\q+4q^{8s+4}\varphi(q^{32s+16})\psi(q^{64s+32})\cdot2q^{8m}\varphi(q^{16s+8})
\psi(q^{64m})
\endalign$$ and so
$$\align&\sum_{n=0}^{\infty}N(2s+1,2s+1,4s+2,8m;8n+4)q^{8n}
\\&=8q^{8s}\psi(q^{16s+8})^2\psi(q^{32s+16})\varphi(q^{32m})
+16q^{8m+8s}\psi(q^{16s+8})^2\psi(q^{32s+16})\psi(q^{64m})
\\&\q+4q^{8s}\varphi(q^{32s+16})\psi(q^{64s+32})\varphi(q^{16s+8})\varphi(q^{32m})
\\&\q+8q^{8s+8m}\varphi(q^{32s+16})\psi(q^{64s+32})\varphi(q^{16s+8})
\psi(q^{64m}).
\endalign$$
Replacing $q$ with $q^{1/8}$ in the above formula we obtain
$$\align&\sum_{n=0}^{\infty}N(2s+1,2s+1,4s+2,8m;8n+4)q^{n}
\\&=8q^{s}\psi(q^{2s+1})^2\psi(q^{4s+2})\varphi(q^{4m})
+16q^{m+s}\psi(q^{2s+1})^2\psi(q^{4s+2})\psi(q^{8m})
\\&\q+4q^{s}\varphi(q^{4s+2})\psi(q^{8s+4})\varphi(q^{2s+1})\varphi(q^{4m})
+8q^{s+m}\varphi(q^{4s+2})\psi(q^{8s+4})\varphi(q^{2s+1})
\psi(q^{8m}).\endalign$$ By (1.8),
$$\varphi(q^{2s+1})\varphi(q^{4s+2})\psi(q^{8s+4})
=\f{\psi(q^{2s+1})^2}{\psi(q^{4s+2})}\cdot
\f{\psi(q^{4s+2})^2}{\psi(q^{8s+4})}\cdot \psi(q^{8s+4})
=\psi(q^{2s+1})^2\psi(q^{4s+2}).$$ Hence
$$\align&\sum_{n=0}^{\infty}N(2s+1,2s+1,4s+2,8m;8n+4)q^{n}
\\&=12q^{s}\psi(q^{2s+1})^2\psi(q^{4s+2})\varphi(q^{4m})
+24q^{s+m}\psi(q^{2s+1})^2\psi(q^{4s+2}) \psi(q^{8m}).\endalign$$ On
the other hand, from (2.1) we have
$$\align&\sum_{n=0}^{\infty}N(2s+1,2s+1,4s+2,8m;2n+1)q^{2n+1}
\\&=4q^{2s+1}\varphi(q^{8s+4})\psi(q^{16s+8})\varphi(q^{4s+2})\varphi(q^{8m}).
\endalign$$
Replacing $q$ with $q^{1/2}$ in the above formula we obtain
$$\align &\sum_{n=0}^{\infty}N(2s+1,2s+1,4s+2,8m;2n+1)q^{n}
\\&=4q^{s}\varphi(q^{4s+2})\psi(q^{8s+4})\varphi(q^{2s+1})
\varphi(q^{4m})
=4q^{s}\psi(q^{2s+1})^2\psi(q^{4s+2})\varphi(q^{4m}).\endalign$$
 Hence,
$$\align&\sum_{n=0}^{\infty}N(2s+1,2s+1,4s+2,8m;8n+4)q^{n}
\\&-3\sum_{n=0}^{\infty}N(2s+1,2s+1,4s+2,8m;2n+1)q^{n}
\\&=24q^{s+m}\psi(q^{2s+1})^2\psi(q^{4s+2})
\psi(q^{8m})
\\&=24q^{s+m}\sum_{n=0}^{\infty}t'(2s+1,2s+1,4s+2,8m;n)q^n
\\&=\f32q^{m+s}
\sum_{n=0}^{\infty}t(2s+1,2s+1,4s+2,8m;n)q^n.\endalign$$
 Comparing
the coefficients of $q^{n+m+s}$ on both sides we obtain the result.

\pro{Theorem 2.2} Let $a\in\{1,3,5,\ldots\}$,
$k,m\in\{0,1,2,\ldots\}$ and $n\in\Bbb N$. Then
$$\align t(a,3a,8k+2,8m+6;n)
&=\f
23N(a,3a,8k+2,8m+6;8n+8k+8m+4a+8)\\&\q-2N(a,3a,8k+2,8m+6;2n+2k+2m+a+2).\endalign$$
\endpro
Proof. Suppose $|q|<1$ and $a=2s+1$. By (1.9),
 $$\aligned &\sum_{n=0}^{\infty}N(2s+1,6s+3,8k+2,8m+6;n)q^{n}
 \\&=\varphi(q^{2s+1})\varphi(q^{6s+3})\varphi(q^{8k+2})
 \varphi(q^{8m+6})
 \\&=\big(\varphi(q^{32s+16})+2q^{8s+4}\psi(q^{64s+32})+2q^{2s+1}\psi(q^{16s+8})\big)
 \\&\q\times\big(\varphi(q^{96s+48})+2q^{24s+12}\psi(q^{192s+96})+2q^{6s+3}\psi(q^{48s+24})\big)
 \\&\q\times\big(\varphi(q^{32k+8})+2q^{8k+2}\psi(q^{64k+16})\big)
 \cdot\big(\varphi(q^{32m+24})+2q^{8m+6}\psi(q^{64m+48})\big)
\\&=
\big(\varphi(q^{32s+16})\varphi(q^{96s+48})+2q^{24s+12}\varphi(q^{32s+16})\psi(q^{192s+96})
+2q^{6s+3}\varphi(q^{32s+16})\psi(q^{48s+24})
\\&\q +2q^{8s+4}\psi(q^{64s+32})\varphi(q^{96s+48})+4q^{32s+16}\psi(q^{64s+32})\psi(q^{192s+96})
\\&\q+4q^{14s+7}\psi(q^{64s+32})\psi(q^{48s+24})
+2q^{2s+1}\psi(q^{16s+8})\varphi(q^{96s+48})
\\&\q+4q^{26s+13}\psi(q^{16s+8})\psi(q^{192s+96})+4q^{8s+4}
\psi(q^{16s+8})\psi(q^{48s+24}) \big)
\\&\q\times \big(\varphi(q^{32k+8})\varphi(q^{32m+24})+
2q^{8m+6}\varphi(q^{32k+8})\psi(q^{64m+48})
+2q^{8k+2}\psi(q^{64k+16})\varphi(q^{32m+24})
\\&\q +4q^{8m+8k+8}\psi(q^{64k+16})\psi(q^{64m+48})\big).
\endaligned\tag 2.2$$
Note that $\varphi(q^{8k_1})^{m_1}\psi(q^{8k_2})^{m_2}
=\sum_{n=0}^{\infty}b_nq^{8n}$ for any nonnegative integers
$k_1,k_2,m_1$ and $m_2$. From (2.2) we deduce that
$$\align&\sum_{n=0}^{\infty}N(2s+1,6s+3,8k+2,8m+6;8n+4)q^{8n+4}
\\&=2q^{24s+12}\varphi(q^{32s+16})\psi(q^{192s+96})
\cdot\varphi(q^{32k+8})\varphi(q^{32m+24})
\\&\q+2q^{24s+12}\varphi(q^{32s+16})\psi(q^{192s+96})
\cdot4q^{8m+8k+8}\psi(q^{64k+16})\psi(q^{64m+48})
\\&\q+2q^{8s+4}\psi(q^{64s+32})\varphi(q^{96s+48})
\cdot\varphi(q^{32k+8})\varphi(q^{32m+24})
\\&\q+2q^{8s+4}\psi(q^{64s+32})\varphi(q^{96s+48})
\cdot4q^{8m+8k+8}\psi(q^{64k+16})\psi(q^{64m+48})
\\&\q+4q^{8s+4}\psi(q^{16s+8})\psi(q^{48s+24})
\cdot\varphi(q^{32k+8})\varphi(q^{32m+24})
\\&\q+4q^{8s+4}\psi(q^{16s+8})\psi(q^{48s+24})
\cdot4q^{8m+8k+8}\psi(q^{64k+16})\psi(q^{64m+48})
\endalign$$ and so
$$\align&\sum_{n=0}^{\infty}N(2s+1,6s+3,8k+2,8m+6;8n+4)q^{8n}
\\&=2q^{24s+8}\varphi(q^{32s+16})\psi(q^{192s+96})\varphi(q^{32k+8})\varphi(q^{32m+24})
\\&\q+8q^{24s+8m+8k+16}\varphi(q^{32s+16})\psi(q^{192s+96})\psi(q^{64k+16})\psi(q^{64m+48})
\\&\q+2q^{8s}\psi(q^{64s+32})\varphi(q^{96s+48})\varphi(q^{32k+8})\varphi(q^{32m+24})
\\&\q+8q^{8m+8k+8s+8}\psi(q^{64s+32})\varphi(q^{96s+48})\psi(q^{64k+16})\psi(q^{64m+48})
\\&\q+4q^{8s}\psi(q^{16s+8})\psi(q^{48s+24})\varphi(q^{32k+8})\varphi(q^{32m+24})
\\&\q+16q^{8m+8k+8s+8}\psi(q^{16s+8})\psi(q^{48s+24})\psi(q^{64k+16})\psi(q^{64m+48}).
\endalign$$
 Replacing $q$ with
$q^{1/8}$ in the above we obtain
$$\align&\sum_{n=0}^{\infty}N(2s+1,6s+3,8k+2,8m+6;8n+4)q^{n}
\\&=2q^{3s+1}\varphi(q^{4s+2})\psi(q^{24s+12})\varphi(q^{4k+1})\varphi(q^{4m+3})
\\&\q+2q^{s}\psi(q^{8s+4})\varphi(q^{12s+6})\varphi(q^{4k+1})\varphi(q^{4m+3})
\\&\q+8q^{3s+m+k+2}\varphi(q^{4s+2})\psi(q^{24s+12})\psi(q^{8k+2})\psi(q^{8m+6})
\\&\q+8q^{m+k+s+1}\psi(q^{8s+4})\varphi(q^{12s+6})\psi(q^{8k+2})\psi(q^{8m+6})
\\&\q+4q^{s}\psi(q^{2s+1})\psi(q^{6s+3})\varphi(q^{4k+1})\varphi(q^{4m+3})
\\&\q+16q^{m+k+s+1}\psi(q^{2s+1})\psi(q^{6s+3})\psi(q^{8k+2})\psi(q^{8m+6}).
\endalign$$
Applying (1.7) we get
$$\aligned&\sum_{n=0}^{\infty}N(2s+1,6s+3,8k+2,8m+6;8n+4)q^{n}
\\&=6q^s\psi(q^{2s+1})\psi(q^{6s+3})\varphi(q^{4k+1})\varphi(q^{4m+3})
+24q^{m+k+s+1}\psi(q^{2s+1})\psi(q^{6s+3})\psi(q^{8k+2})\psi(q^{8m+6}).
\endaligned$$
By (1.6),
$$\align&\sum_{n=0}^{\infty}N(2s+1,6s+3,8k+2,8m+6;n)q^{n}
\\&=\varphi(q^{2s+1})\varphi(q^{6s+3})\varphi(q^{8k+2})\varphi(q^{8m+6})
\\&=\big(\varphi(q^{8s+4})+2q^{2s+1}\psi(q^{16s+8})\big)\big(\varphi(q^{24s+12})+2q^{6s+3}
\psi(q^{48s+24})\big)\varphi(q^{8k+2})\varphi(q^{8m+6})
\endalign$$ and so
$$\align&\sum_{n=0}^{\infty}N(2s+1,6s+3,8k+2,8m+6;2n+1)q^{2n+1}
\\&=2q^{6s+3}\varphi(q^{8s+4})\psi(q^{48s+24})\varphi(q^{8k+2})\varphi(q^{8m+6})
\\&\q+2q^{2s+1}\psi(q^{16s+8})\varphi(q^{24s+12})\varphi(q^{8k+2})\varphi(q^{8m+6}).
\endalign$$ Replacing $q$ with $q^{1/2}$ in the above formula we obtain
$$\align&\sum_{n=0}^{\infty}N(2s+1,6s+3,8k+2,8m+6;2n+1)q^{n}
\\&=2q^{3s+1}\varphi(q^{4s+2})\psi(q^{24s+12})\varphi(q^{4k+1})\varphi(q^{4m+3})
+2q^{s}\psi(q^{8s+4})\varphi(q^{12s+6})\varphi(q^{4k+1})\varphi(q^{4m+3}).
\endalign$$ Now applying (1.7) we get
$$\sum_{n=0}^{\infty}N(2s+1,6s+3,8k+2,8m+6;2n+1)q^{n}
=2q^s\psi(q^{2s+1})\psi(q^{6s+3})\varphi(q^{4k+1})\varphi(q^{4m+3}).$$
Thus,
$$\align&\sum_{n=0}^{\infty}N(2s+1,6s+3,8k+2,8m+6;8n+4)q^{n}
\\&\q-3\sum_{n=0}^{\infty}N(2s+1,6s+3,8k+2,8m+6;2n+1)q^{n}
\\&=24q^{m+k+s+1}\psi(q^{2s+1})\psi(q^{6s+3})\psi(q^{8k+2})\psi(q^{8m+6})
\\&=24q^{m+k+s+1}\sum_{n=0}^{\infty}t'(2s+1,6s+3,8k+2,8m+6;n)q^{n}
\\&=\f32q^{m+k+s+1}\sum_{n=0}^{\infty}
t(2s+1,6s+3,8k+2,8m+6;n)q^{n}.\endalign$$
 Comparing the
coefficients of $q^{m+n+k+s+1}$ on both sides we obtain the result.

\pro{Theorem 2.3} Let $a\in\{1,3,5,\ldots\}$, $m\in\{0,1,2,\ldots\}$
and $n\in\Bbb N$. If $n\e m+\f{a-1}2\mod 2$, then
$$\align t(a,3a,8m+4,8m+4;n)
&=\f
23N(a,3a,8m+4,8m+4;8n+16m+4a+8)\\&\q-2N(a,3a,8m+4,8m+4;2n+4m+a+2).\endalign$$
\endpro
Proof. Suppose $|q|<1$ and $a=2s+1$. Using (1.9) we see that
 $$\aligned&\sum_{n=0}^{\infty}N(2s+1,6s+3,8m+4,8m+4;n)q^{n}
 \\&=\varphi(q^{2s+1})\varphi(q^{6s+3})\varphi(q^{8m+4})^2
 \\&=\big(\varphi(q^{32s+16})+2q^{8s+4}\psi(q^{64s+32})
 +2q^{2s+1}\psi(q^{16s+8})\big)
 \\&\q\times\big(\varphi(q^{96s+48})+2q^{24s+12}
 \psi(q^{192s+96})
 +2q^{6s+3}\psi(q^{48s+24})\big)
 \\&\q\times\big(\varphi(q^{32m+16})+2q^{8m+4}
 \psi(q^{64m+32})\big)^2
\\&=
\big(\varphi(q^{32s+16})\varphi(q^{96s+48})+2q^{24s+12}
\varphi(q^{32s+16})\psi(q^{192s+96})
+2q^{6s+3}\varphi(q^{32s+16})\psi(q^{48s+24})
\\&\q +2q^{8s+4}\psi(q^{64s+32})\varphi(q^{96s+48})+4q^{32s+16}\psi(q^{64s+32})\psi(q^{192s+96})
\\&\q+4q^{14s+7}\psi(q^{64s+32})\psi(q^{48s+24})
+2q^{2s+1}\psi(q^{16s+8})\varphi(q^{96s+48})
\\&\q+4q^{26s+13}\psi(q^{16s+8})\psi(q^{192s+96})+4q^{8s+4}
\psi(q^{16s+8})\psi(q^{48s+24}) \big)
\\&\q\times \big(\varphi(q^{32m+16})^2+
4q^{8m+4}\varphi(q^{32m+16})\psi(q^{64m+32})
+4q^{16m+8}\psi(q^{64m+32})^2 \big).
\endaligned\tag 2.3$$
Note that $\varphi(q^{8k_1})^{m_1}\psi(q^{8k_2})^{m_2}
=\sum_{n=0}^{\infty}b_nq^{8n}$ for any nonnegative integers
$k_1,k_2,m_1$ and $m_2$. From (2.3) we deduce that
$$\align&\sum_{n=0}^{\infty}N(2s+1,6s+3,8m+4,8m+4;8n+4)q^{8n+4}
\\&=\varphi(q^{32s+16})\varphi(q^{96s+48})
\cdot4q^{8m+4}\varphi(q^{32m+16})\psi(q^{64m+32})
\\&\q+2q^{24s+12}\varphi(q^{32s+16})\psi(q^{192s+96})
(\varphi(q^{32m+16})^2 +4q^{16m+8}\psi(q^{64m+32})^2)
\\&\q+2q^{8s+4}\psi(q^{64s+32})\varphi(q^{96s+48})
(\varphi(q^{32m+16})^2 +4q^{16m+8}\psi(q^{64m+32})^2)
\\&\q+4q^{32s+16}\psi(q^{64s+32})\psi(q^{192s+96})
\cdot4q^{8m+4}\varphi(q^{32m+16})\psi(q^{64m+32})
\\&\q+4q^{8s+4}\psi(q^{16s+8})\psi(q^{48s+24})
(\varphi(q^{32m+16})^2
 +4q^{16m+8}\psi(q^{64m+32})^2) \endalign$$ and
so
$$\align&\sum_{n=0}^{\infty}N(2s+1,6s+3,8m+4,8m+4;8n+4)q^{8n}
\\&=4q^{8m}\varphi(q^{32s+16})\varphi(q^{96s+48})
\varphi(q^{32m+16})\psi(q^{64m+32})
\\&\q+2q^{24s+8}\varphi(q^{32s+16})\psi(q^{192s+96})
\varphi(q^{32m+16})^2
\\&\q+8q^{24s+16m+16}\varphi(q^{32s+16})\psi(q^{192s+96})
\psi(q^{64m+32})^2
\\&\q+2q^{8s}\psi(q^{64s+32})\varphi(q^{96s+48})
\varphi(q^{32m+16})^2
\\&\q+8q^{8s+16m+8}\psi(q^{64s+32})\varphi(q^{96s+48})
\psi(q^{64m+32})^2
\\&\q+16q^{32s+8m+16}\psi(q^{64s+32})\psi(q^{192s+96})
\varphi(q^{32m+16})\psi(q^{64m+32})
\\&\q+4q^{8s}\psi(q^{16s+8})\psi(q^{48s+24})
\varphi(q^{32m+16})^2
\\&\q+16q^{8s+16m+8}\psi(q^{16s+8})\psi(q^{48s+24})
\psi(q^{64m+32})^2 .\endalign$$
 Replacing $q$ with
$q^{1/8}$  we then obtain
$$\align&\sum_{n=0}^{\infty}N(2s+1,6s+3,8m+4,8m+4;8n+4)q^{n}
\\&=4q^{m}\varphi(q^{4s+2})\varphi(q^{12s+6})
\varphi(q^{4m+2})\psi(q^{8m+4})
\\&\q+2q^{3s+1}\varphi(q^{4s+2})\psi(q^{24s+12}) \varphi(q^{4m+2})^2
+2q^{s}\psi(q^{8s+4})\varphi(q^{12s+6}) \varphi(q^{4m+2})^2
\\&\q+8q^{3s+2m+2}\varphi(q^{4s+2})\psi(q^{24s+12}) \psi(q^{8m+4})^2
+8q^{s+2m+1}\psi(q^{8s+4})\varphi(q^{12s+6}) \psi(q^{8m+4})^2
\\&\q+16q^{4s+m+2}\psi(q^{8s+4})\psi(q^{24s+12})
\varphi(q^{4m+2})\psi(q^{8m+4})
\\&\q+4q^{s}\psi(q^{2s+1})\psi(q^{6s+3})
\varphi(q^{4m+2})^2 +16q^{s+2m+1}\psi(q^{2s+1})\psi(q^{6s+3})
\psi(q^{8m+4})^2 .\endalign$$ By (1.7),
$$\psi(q^{8s+4})\varphi(q^{12s+6})+q^{2s+1}\varphi(q^{4s+2})
\psi(q^{24s+12})=\psi(q^{2s+1})\psi(q^{6s+3}).\tag 2.4$$
 Thus,
$$\aligned&\sum_{n=0}^{\infty}N(2s+1,6s+3,8m+4,8m+4;8n+4)q^{n}
\\&=4q^{m}\varphi(q^{4s+2})\varphi(q^{12s+6})
\varphi(q^{4m+2})\psi(q^{8m+4})
\\&\q+16q^{4s+m+2}\psi(q^{8s+4})\psi(q^{24s+12})
\varphi(q^{4m+2})\psi(q^{8m+4})
\\&\q+6q^{s}\psi(q^{2s+1})\psi(q^{6s+3})
\varphi(q^{4m+2})^2 +24q^{s+2m+1}\psi(q^{2s+1})\psi(q^{6s+3})
\psi(q^{8m+4})^2 .\endaligned\tag 2.5$$ On the other hand, using
(1.9) we see that
$$\align&\sum_{n=0}^{\infty}N(2s+1,6s+3,8m+4,8m+4;n)q^{n}
\\&=\varphi(q^{2s+1})\varphi(q^{6s+3})\varphi(q^{8m+4})^2
\\&=\big(\varphi(q^{8s+4})+2q^{2s+1}\psi(q^{16s+8})\big)
\big(\varphi(q^{24s+12})
+2q^{6s+3}\psi(q^{48s+24})\big)\varphi(q^{8m+4})^2
\endalign$$ and so
$$\align&\sum_{n=0}^{\infty}N(2s+1,6s+3,8m+4,8m+4;2n+1)q^{2n+1}
\\&=\big(2q^{6s+3}\varphi(q^{8s+4})\psi(q^{48s+24})
+2q^{2s+1}\psi(q^{16s+8})\varphi(q^{24s+12})\big)\varphi(q^{8m+4})^2.
\endalign$$ Replacing $q$ with $q^{1/2}$ we then obtain
$$\align&\sum_{n=0}^{\infty}N(2s+1,6s+3,8m+4,8m+4;2n+1)q^{n}
\\&=\big(2q^{3s+1}\varphi(q^{4s+2})\psi(q^{24s+12})
+2q^{s}\psi(q^{8s+4})\varphi(q^{12s+6})\big)\varphi(q^{4m+2})^2.
\endalign$$ Now applying (2.4) we get
$$\sum_{n=0}^{\infty}N(2s+1,6s+3,8m+4,8m+4;2n+1)q^{n}
=2q^{s}\psi(q^{2s+1})\psi(q^{6s+3})\varphi(q^{4m+2})^2.\tag 2.6$$
From (2.5) and (2.6) we deduce that
$$\align&\sum_{n=0}^{\infty}N(2s+1,6s+3,8m+4,8m+4;8n+4)q^{n}
\\&-3\sum_{n=0}^{\infty}N(2s+1,6s+3,8m+4,8m+4;2n+1)q^{n}
\\&=4q^{m}\varphi(q^{4s+2})\varphi(q^{12s+6})
\varphi(q^{4m+2})\psi(q^{8m+4})
\\&\q+16q^{4s+m+2}\psi(q^{8s+4})\psi(q^{24s+12})
\varphi(q^{4m+2})\psi(q^{8m+4})
\\&\q +24q^{s+2m+1}\psi(q^{2s+1})\psi(q^{6s+3})
\psi(q^{8m+4})^2
\\&=4q^{m}\varphi(q^{4s+2})\varphi(q^{12s+6})
\varphi(q^{4m+2})\psi(q^{8m+4})
\\&\q+16q^{4s+m+2}\psi(q^{8s+4})\psi(q^{24s+12})
\varphi(q^{4m+2})\psi(q^{8m+4})
\\&\q+24q^{s+2m+1}\sum_{n=0}^{\infty}t'(2s+1,6s+3,8m+4,8m+4;n)q^n.
\endalign$$
 Suppose $n\e m+s\mod 2$. Then $s+2m+1+n\e m+1\mod 2$. Comparing the coefficients
of $q^{s+2m+1+n}$ in the above expansion we obtain
$$\align&N(2s+1,6s+3,8m+4,8m+4;8(2m+n+s+1)+4)
\\&-3N(2s+1,6s+3,8m+4,8m+4;2(2m+n+s+1)+1)
\\&=24t'(2s+1,6s+3,8m+4,8m+4;n)=\f32t(2s+1,6s+3,8m+4,8m+4;n).\endalign$$ This completes the
proof.

\pro{Theorem 2.4} Let $a\in\{1,3,5,\ldots\}$,
$k,m\in\{0,1,2,\ldots\}$ and $n\in\Bbb N$. If $n\e \f{a-1}2\mod 2$,
then

$$\align t(a,3a,16k+4,16m+4;n)
&=\f
23N(a,3a,16k+4,16m+4;8n+16k+16m+4a+8)\\&\q-2N(a,3a,16k+4,16m+4;2n+4k+4m+a+2).\endalign$$
\endpro
Proof. Suppose $|q|<1$ and $a=2s+1$. Using (1.9) we see that
 $$\aligned&\sum_{n=0}^{\infty}N(2s+1,6s+3,16k+4,16m+4;n)q^{n}
 \\&=\varphi(q^{2s+1})\varphi(q^{6s+3})\varphi(q^{16k+4})
 \varphi(q^{16m+4})
\\&=\big(\varphi(q^{32s+16})+2q^{8s+4}\psi(q^{64s+32})+2q^{2s+1}\psi(q^{16s+8})\big)
\\&\q\times \big(\varphi(q^{96s+48})+2q^{24s+12}\psi(q^{192s+96})+2q^{6s+3}\psi(q^{48s+24})\big)
\\&\q\times
\big(\varphi(q^{64k+16})+2q^{16k+4}\psi(q^{128k+32})\big)\cdot
\big(\varphi(q^{64m+16})+2q^{16m+4}\psi(q^{128m+32})\big)
\\&=
\big(\varphi(q^{32s+16})\varphi(q^{96s+48})+2q^{24s+12}
\varphi(q^{32s+16})\psi(q^{192s+96})
+2q^{6s+3}\varphi(q^{32s+16})\psi(q^{48s+24})
\\&\q +2q^{8s+4}\psi(q^{64s+32})\varphi(q^{96s+48})+4q^{32s+16}\psi(q^{64s+32})\psi(q^{192s+96})
\\&\q+4q^{14s+7}\psi(q^{64s+32})\psi(q^{48s+24})
+2q^{2s+1}\psi(q^{16s+8})\varphi(q^{96s+48})
\\&\q+4q^{26s+13}\psi(q^{16s+8})\psi(q^{192s+96})+4q^{8s+4}
\psi(q^{16s+8})\psi(q^{48s+24}) \big)
\\&\q\times \big(\varphi(q^{64k+16})\varphi(q^{64m+16})+
2q^{16m+4}\varphi(q^{64k+16})\psi(q^{128m+32})
\\&\q+2q^{16k+4}\psi(q^{128k+32})\varphi(q^{64m+16}) +
4q^{16k+16m+8}\psi(q^{128k+32})\psi(q^{128m+32})\big).
\endaligned\tag 2.7$$
Note that $\varphi(q^{8k_1})^{m_1}\psi(q^{8k_2})^{m_2}
=\sum_{n=0}^{\infty}b_nq^{8n}$ for any nonnegative integers
$k_1,k_2,m_1$ and $m_2$. From (2.7) we deduce that
$$\align&\sum_{n=0}^{\infty}N(2s+1,6s+3,16k+4,16m+4;8n+4)q^{8n+4}
\\&=\varphi(q^{32s+16})\varphi(q^{96s+48})
\big(2q^{16m+4}\varphi(q^{64k+16})\psi(q^{128m+32}) +
2q^{16k+4}\psi(q^{128k+32})\varphi(q^{64m+16})\big)
\\&\q+2q^{24s+12}\varphi(q^{32s+16})\psi(q^{192s+96})
\big(\varphi(q^{64k+16})\varphi(q^{64m+16})
\\&\qq+4q^{16k+16m+8}\psi(q^{128k+32})\psi(q^{128m+32})\big)
\\&\q+2q^{8s+4}\psi(q^{64s+32})\varphi(q^{96s+48})
\big(\varphi(q^{64k+16})\varphi(q^{64m+16})
\\&\qq+4q^{16k+16m+8}\psi(q^{128k+32})\psi(q^{128m+32})\big)
\\&\q+4q^{32s+16}\psi(q^{64s+32})\psi(q^{192s+96})
\big(2q^{16m+4}\varphi(q^{64k+16})\psi(q^{128m+32})
\\&\qq+2q^{16k+4}\psi(q^{128k+32})\varphi(q^{64m+16})\big)
\\&\q+4q^{8s+4}\psi(q^{16s+8})\psi(q^{48s+24})
\big(\varphi(q^{64k+16})\varphi(q^{64m+16})
\\&\qq+4q^{16k+16m+8}\psi(q^{128k+32})\psi(q^{128m+32})\big)\endalign$$
and so
$$\align&\sum_{n=0}^{\infty}N(2s+1,6s+3,16k+4,16m+4;8n+4)q^{8n}
\\&=2q^{16m}\varphi(q^{32s+16})\varphi(q^{96s+48})
\varphi(q^{64k+16})\psi(q^{128m+32})
\\&\q+2q^{16k}\varphi(q^{32s+16})\varphi(q^{96s+48})
\psi(q^{128k+32})\varphi(q^{64m+16})
\\&\q+2q^{24s+8}\varphi(q^{32s+16})\psi(q^{192s+96})
\varphi(q^{64k+16})\varphi(q^{64m+16})
\\&\q+8q^{24s+16k+16m+16}\varphi(q^{32s+16})\psi(q^{192s+96})
\psi(q^{128k+32})\psi(q^{128m+32})
\\&\q+2q^{8s}\psi(q^{64s+32})\varphi(q^{96s+48})
\varphi(q^{64k+16})\varphi(q^{64m+16})
\\&\q+8q^{8s+16k+16m+8}\psi(q^{64s+32})\varphi(q^{96s+48})
\psi(q^{128k+32})\psi(q^{128m+32})
\\&\q+8q^{32s+16m+16}\psi(q^{64s+32})\psi(q^{192s+96})
\varphi(q^{64k+16})\psi(q^{128m+32})
\\&\q+8q^{32s+16k+16}\psi(q^{64s+32})\psi(q^{192s+96})
\psi(q^{128k+32})\varphi(q^{64m+16})
\\&\q+4q^{8s}\psi(q^{16s+8})\psi(q^{48s+24})
\varphi(q^{64k+16})\varphi(q^{64m+16})
\\&\q+16q^{8s+16k+16m+8}\psi(q^{16s+8})\psi(q^{48s+24})
\psi(q^{128k+32})\psi(q^{128m+32}).\endalign$$
 Replacing $q$ with
$q^{1/8}$ in the above we obtain
$$\align&\sum_{n=0}^{\infty}N(2s+1,6s+3,16k+4,16m+4;8n+4)q^{n}
\\&=2q^{2m}\varphi(q^{4s+2})\varphi(q^{12s+6})
\varphi(q^{8k+2})\psi(q^{16m+4})
\\&\q+2q^{2k}\varphi(q^{4s+2})\varphi(q^{12s+6})
\psi(q^{16k+4})\varphi(q^{8m+2})
\\&\q+2q^{3s+1}\varphi(q^{4s+2})\psi(q^{24s+12})
\varphi(q^{8k+2})\varphi(q^{8m+2})
\\&\q+2q^{s}\psi(q^{8s+4})\varphi(q^{12s+6})
\varphi(q^{8k+2})\varphi(q^{8m+2})
\\&\q+8q^{3s+2k+2m+2}\varphi(q^{4s+2})\psi(q^{24s+12})
\psi(q^{16k+4})\psi(q^{16m+4})
\\&\q+8q^{s+2k+2m+1}\psi(q^{8s+4})\varphi(q^{12s+6})
\psi(q^{16k+4})\psi(q^{16m+4})
\\&\q+8q^{4s+2m+2}\psi(q^{8s+4})\psi(q^{24s+12})
\varphi(q^{8k+2})\psi(q^{16m+4})
\\&\q+8q^{4s+2k+2}\psi(q^{8s+4})\psi(q^{24s+12})
\psi(q^{16k+4})\varphi(q^{8m+2})
\\&\q+4q^{s}\psi(q^{2s+1})\psi(q^{6s+3})
\varphi(q^{8k+2})\varphi(q^{8m+2})
\\&\q+16q^{s+2k+2m+1}\psi(q^{2s+1})\psi(q^{6s+3})
\psi(q^{16k+4})\psi(q^{16m+4}).\endalign$$
 By (1.7),
$$\psi(q^{8s+4})\varphi(q^{12s+6})+q^{2s+1}\varphi(q^{4s+2})
\psi(q^{24s+12})=\psi(q^{2s+1})\psi(q^{6s+3}).$$
 Thus,
$$\aligned&\sum_{n=0}^{\infty}N(2s+1,6s+3,16k+4,16m+4;8n+4)q^{n}
\\&=2q^{2m}\varphi(q^{4s+2})\varphi(q^{12s+6})
\varphi(q^{8k+2})\psi(q^{16m+4})
\\&\q+2q^{2k}\varphi(q^{4s+2})\varphi(q^{12s+6})
\psi(q^{16k+4})\varphi(q^{8m+2})
\\&\q+8q^{4s+2m+2}\psi(q^{8s+4})\psi(q^{24s+12})
\varphi(q^{8k+2})\psi(q^{16m+4})
\\&\q+8q^{4s+2k+2}\psi(q^{8s+4})\psi(q^{24s+12})
\psi(q^{16k+4})\varphi(q^{8m+2})
\\&\q+6q^{s}\psi(q^{2s+1})\psi(q^{6s+3})
\varphi(q^{8k+2})\varphi(q^{8m+2})
\\&\q+24q^{s+2k+2m+1}\psi(q^{2s+1})\psi(q^{6s+3})
\psi(q^{16k+4})\psi(q^{16m+4}).\endaligned\tag 2.8$$ On the other
hand, using (1.9) we see that
$$\align&\sum_{n=0}^{\infty}N(2s+1,6s+3,16k+4,16m+4;n)q^{n}
\\&=\varphi(q^{2s+1})\varphi(q^{6s+3})\varphi(q^{16k+4})\varphi(q^{16m+4})
\\&=\Big(\varphi(q^{8s+4})+2q^{2s+1}\psi(q^{16s+8})\Big)\Big(\varphi(q^{24s+12})
+2q^{6s+3}\psi(q^{48s+24})\Big)\varphi(q^{16k+4})\varphi(q^{16m+4})
\endalign$$
and so
$$\align&\sum_{n=0}^{\infty}N(2s+1,6s+3,16k+4,16m+4;2n+1)q^{2n+1}
\\&=2q^{6s+3}\varphi(q^{8s+4})\psi(q^{48s+24})\varphi(q^{16k+4})
\varphi(q^{16m+4})
\\&\q+2q^{2s+1}\psi(q^{16s+8})\varphi(q^{24s+12})\varphi(q^{16k+4})
\varphi(q^{16m+4}).
\endalign$$
Replacing $q$ with $q^{1/2}$ in the above formula we obtain
$$\align&\sum_{n=0}^{\infty}N(2s+1,6s+3,16k+4,16m+4;2n+1)q^{n}
\\&=2q^{3s+1}\varphi(q^{4s+2})\psi(q^{24s+12})\varphi(q^{8k+2})
\varphi(q^{8m+2})
+2q^{s}\psi(q^{8s+4})\varphi(q^{12s+6})\varphi(q^{8k+2})
\varphi(q^{8m+2}).
\endalign$$
Now applying (1.7) we get
$$\aligned&\sum_{n=0}^{\infty}N(2s+1,6s+3,16k+4,16m+4;2n+1)q^{n}
\\&=2q^{s}\psi(q^{2s+1})\psi(q^{6s+3})\varphi(q^{12s+6})
\varphi(q^{8k+2}).\endaligned\tag 2.9$$
From (2.8) and (2.9) we
deduce that
$$\aligned&\sum_{n=0}^{\infty}N(2s+1,6s+3,16k+4,16m+4;8n+4)q^{n}
\\&\q-3\sum_{n=0}^{\infty}N(2s+1,6s+3,16k+4,16m+4;2n+1)q^{n}
\\&=2q^{2m}\varphi(q^{4s+2})\varphi(q^{12s+6})
\varphi(q^{8k+2})\psi(q^{16m+4})
\\&\q+2q^{2k}\varphi(q^{4s+2})\varphi(q^{12s+6})
\psi(q^{16k+4})\varphi(q^{8m+2})
\\&\q+8q^{4s+2m+2}\psi(q^{8s+4})\psi(q^{24s+12})
\varphi(q^{8k+2})\psi(q^{16m+4})
\\&\q+8q^{4s+2k+2}\psi(q^{8s+4})\psi(q^{24s+12})
\psi(q^{16k+4})\varphi(q^{8m+2})
\\&\q+24q^{s+2k+2m+1}\psi(q^{2s+1})\psi(q^{6s+3})
\psi(q^{16k+4})\psi(q^{16m+4}).\endaligned$$ Suppose $n\e s\mod 2$.
Then $n+s+2k+2m+1\e 1\mod 2$. Comparing the coefficients of
$q^{n+s+2k+2m+1}$ in the above expansion we obtain
$$\align&N(2s+1,6s+3,16k+4,16m+4;8(n+s+2k+2m+1)+4)
\\&-3N(2s+1,6s+3,16k+4,16m+4;2(n+s+2k+2m+1)+1)
\\&=24t'(2s+1,6s+3,16k+4,16m+4;n)=\f32t(2s+1,6s+3,16k+4,16m+4;n).\endalign$$ This completes the
proof.

\section*{3. Formulas for $t(1,1,2,8;n)$, $t(1,1,2,16;n)$, $t(1,2,3,6;n)$
and $t(1,3,4,12;n)$}
 \pro{Lemma 3.1 ([AALW2, Theorem 4.3])} Suppose
$n\in\Bbb N$, $n=2^{\a}n_1$ and $2\nmid n_1$. Then
$$\aligned N(1,1,2,8;n)=\cases 2\sigma(n_1)+2\sls 2{n_1}
\sum\limits\Sb (r,s)\in\Bbb Z^2, 4\mid r-1\\n_1=r^2+4s^2\endSb r
&\t{if $n\e1\mod 4$,}
\\2\sigma(n_1)&\t{if $n\e3\mod 4$,}
\\12\sigma(n_1)&\t{if $n\e 4\mod 8$.}
\endcases\endaligned$$
\endpro
\pro{Theorem 3.1} Let $n\in\Bbb N$. Then
$$\aligned t(1,1,2,8;n)=\cases 4\sigma(2n+3)
&\t{if $2\mid n$,} \\4\sigma(2n+3)+4(-1)^{\f {n-1}2}\sum\limits\Sb
(r,s)\in\Bbb Z^2, 4\mid r-1\\2n+3=r^2+4s^2\endSb r &\t{if $2\nmid
n$.}\endcases\endaligned$$\endpro Proof. Taking $a=m=1$ in
 Theorem 2.1 we see that
 $$t(1,1,2,8;n)=\f 23N(1,1,2,8;8n+12)-2N(1,1,2,8;2n+3).$$
 Now applying Lemma 3.1 we deduce the result.
 \pro{Lemma 3.2 ([AALW2, Theorem 4.15])} Let $n\in\Bbb N$ and
$n=2^{\a}n_1$ with $2\nmid n_1$. Then
$$\aligned N(1,1,2,16;n)=\cases 2\sum_{d\mid n_1} \f {n_1}d\Ls 2d+2\sum\limits\Sb
(r,s)\in\Bbb Z^2, 4\mid r-1\\n_1=r^2+2s^2\endSb r &\t{if $n\e 1\mod
2$,}\\12\sum_{d\mid n_1} \f {n_1}d\Ls 2d&\t{if $n\e 4\mod
8$.}\endcases\endaligned$$
\endpro
\pro{Theorem 3.2} Suppose $n\in\Bbb N$. Then
$$t(1,1,2,16;n)=4\sum_{d\mid {2n+5}} \f {2n+5}d\Ls 2d-4\sum\limits\Sb
(r,s)\in\Bbb Z^2, 4\mid r-1\\2n+5=r^2+2s^2\endSb r .$$\endpro

Proof. Taking $a=1$ and $m=2$ in
 Theorem 2.1 we see that
 $$t(1,1,2,16;n)=\f 23N(1,1,2,16;8n+20)-2N(1,1,2,16;2n+5).$$
 Now applying Lemma 3.2 we deduce the result.

\pro{Lemma 3.3 ([AALW1, Theorem 1.15])} Let $n\in\Bbb N$,
$n=2^{\a}3^{\beta}n_1$ and $\t gcd (n_1,6)=1.$ Then
$$\aligned N(1,2,3,6;n)=\cases (3^{\beta+1}-2)\sigma(n_1)+a(n)
&\t{if $n\e1\mod 2$,}
\\6(3^{\beta+1}-2)\sigma(n_1)&\t{if $n\e0\mod 4$.}
\endcases\endaligned$$
\endpro
\pro{Theorem 3.3} Suppose $n\in\Bbb N$ and $2n+3=3^{\beta}n_1$ with
$n_1\in\Bbb N$ and $3\nmid n_1$. Then
$$t(1,2,3,6;n)=2(3^{\beta+1}-2)\sigma(n_1)-2a(2n+3).$$\endpro
Proof. Taking $a=1$ and $k=m=0$ in
 Theorem 2.2 we see that
 $$t(1,2,3,6;n)=\f 23N(1,2,3,6;8n+12)-2N(1,2,3,6;2n+3).$$
 Now applying Lemma 3.3 we deduce the result.

  \pro{Lemma 3.4 ([AALW1, Theorem 1.17])} Let
$n\in\Bbb N$ and $n=2^{\a}3^{\beta}n_1$ with $gcd(n_1,6)=1$. Then
$$\aligned N(1,3,4,12;n)=\cases 8\sigma(n_1)
&\t{if $n\e4\mod8$,}
\\\sigma(n_1)+a(n)&\t{if $n\e1\mod4$,}
\\\sigma(n_1)-a(n)&\t{if $n\e3\mod4$.}\endcases\endaligned$$
\endpro
\pro{Theorem 3.4} Suppose $n\in\Bbb N$ and $2n+5=3^{\beta}n_1$ with
$3\nmid n_1$. Then
$$t(1,3,4,12;n)= 2(\sigma(n_1)-(-1)^na(2n+5)).$$
\endpro Proof. Suppose $|q|<1$.
Then clearly
 $$\sum_{n=0}^{\infty}N(1,3,4,12;n)q^{n}
 =\varphi(q)\varphi(q^3)\varphi(q^4)\varphi(q^{12}).$$
By (1.9),
$$\align&\varphi(q)\varphi(q^3)\varphi(q^4)\varphi(q^{12})
\\&=\big(\varphi(q^{16})+2q^4\psi(q^{32})+2q\psi(q^{8})\big)
\\&\q\times\big(\varphi(q^{48})+2q^{12}\psi(q^{96})+2q^3\psi(q^{24})\big)
\\&\q\times\big(\varphi(q^{16})+2q^4\psi(q^{32})\big)\cdot
(\varphi(q^{48})+2q^{12}\psi(q^{96})\big)
\\&=\big(\varphi(q^{16})\varphi(q^{48})+2q^{12}\varphi(q^{16})\psi(q^{96})
+2q^3\varphi(q^{16})\psi(q^{24})
\\&\q+2q^4\psi(q^{32})\varphi(q^{48})+4q^{16}\psi(q^{32})\psi(q^{96})
+4q^7\psi(q^{32})\psi(q^{24})
\\&\q+2q\psi(q^{8})\varphi(q^{48})+4q^{13}\psi(q^{8})\psi(q^{96})
+4q^4\psi(q^{8})\psi(q^{24})\big)
\\&\q\times
\big(\varphi(q^{16})\varphi(q^{48})+2q^{12}\varphi(q^{16})
\psi(q^{96})+2q^{4}\psi(q^{32})\varphi(q^{48})
\\&\q+4q^{16}\psi(q^{96})\psi(q^{32})\big).
\endalign$$
Note that $\varphi(q^{8k_1})^{m_1}\psi(q^{8k_2})^{m_2}
=\sum_{n=0}^{\infty}b_nq^{8n}$ for $|q|<1$ and any nonnegative
integers $k_1,k_2,m_1$ and $m_2$. From the above we deduce that
$$\align&\sum_{n=0}^{\infty}N(1,3,4,12;8n+4)q^{8n+4}
\\&=\varphi(q^{16})\varphi(q^{48})\cdot2q^{12}\varphi(q^{16})\psi(q^{96})
+\varphi(q^{16})\varphi(q^{48})\cdot2q^{4}\psi(q^{32})\varphi(q^{48})
\\&\q+2q^{12}\varphi(q^{16})\psi(q^{96})\cdot\varphi(q^{16})\varphi(q^{48})
+2q^{12}\varphi(q^{16})\psi(q^{96})\cdot4q^{16}\psi(q^{96})\psi(q^{32})
\\&\q+2q^4\psi(q^{32})\varphi(q^{48})\cdot\varphi(q^{16})\varphi(q^{48})
+2q^4\psi(q^{32})\varphi(q^{48})\cdot4q^{16}\psi(q^{96})\psi(q^{32})
\\&\q+4q^{16}\psi(q^{32})\psi(q^{96})\cdot2q^{12}\varphi(q^{16})\psi(q^{96})
+4q^{16}\psi(q^{32})\psi(q^{96})\cdot2q^{4}\psi(q^{32})\varphi(q^{48})
\\&\q+4q^4\psi(q^{8})\psi(q^{24})\cdot\varphi(q^{16})\varphi(q^{48})
+4q^4\psi(q^{8})\psi(q^{24})\cdot4q^{16}\psi(q^{96})\psi(q^{32})
\endalign$$ and so

$$\align&\sum_{n=0}^{\infty}N(1,3,4,12;8n+4)q^{8n}
\\&=2q^{8}\varphi(q^{16})\varphi(q^{48})\varphi(q^{16})\psi(q^{96})
+2\varphi(q^{16})\varphi(q^{48})\psi(q^{32})\varphi(q^{48})
\\&\q+2q^{8}\varphi(q^{16})\psi(q^{96})\varphi(q^{16})\varphi(q^{48})
+8q^{24}\varphi(q^{16})\psi(q^{96})\psi(q^{96})\psi(q^{32})
\\&\q+2\psi(q^{32})\varphi(q^{48})\varphi(q^{16})\varphi(q^{48})
+8q^{16}\psi(q^{32})\varphi(q^{48})\psi(q^{96})\psi(q^{32})
\\&\q+8q^{24}\psi(q^{32})\psi(q^{96})\varphi(q^{16})\psi(q^{96})
+8q^{16}\psi(q^{32})\psi(q^{96})\psi(q^{32})\varphi(q^{48})
\\&\q+4\psi(q^{8})\psi(q^{24})\varphi(q^{16})\varphi(q^{48})
+16q^{16}\psi(q^{8})\psi(q^{24})\psi(q^{96})\psi(q^{32}).
\endalign$$
Replacing $q$ with $q^{1/8}$ in the above we obtain
$$\align&\sum_{n=0}^{\infty}N(1,3,4,12;8n+4)q^{n}
\\&=2q\varphi(q^{2})\varphi(q^{6})\varphi(q^{2})\psi(q^{12})
+2\varphi(q^{2})\varphi(q^{6})\psi(q^{4})\varphi(q^{6})
\\&\q+2q\varphi(q^{2})\psi(q^{12})\varphi(q^{2})\varphi(q^{6})
+8q^{3}\varphi(q^{2})\psi(q^{12})\psi(q^{12})\psi(q^{4})
\\&\q+2\psi(q^{4})\varphi(q^{6})\varphi(q^{2})\varphi(q^{6})
+8q^{2}\psi(q^{4})\varphi(q^{6})\psi(q^{12})\psi(q^{4})
\\&\q+8q^{3}\psi(q^{4})\psi(q^{12})\varphi(q^{2})\psi(q^{12})
+8q^{2}\psi(q^{4})\psi(q^{12})\psi(q^{4})\varphi(q^{6})
\\&\q+4\psi(q)\psi(q^{3})\varphi(q^{2})\varphi(q^{6})
+16q^{2}\psi(q)\psi(q^{3})\psi(q^{12})\psi(q^{4})
\\&=4q\varphi(q^{2})\varphi(q^{6})\varphi(q^{2})\psi(q^{12})
+4\varphi(q^{2})\varphi(q^{6})\psi(q^{4})\varphi(q^{6})
\\&\q+4\psi(q)\psi(q^{3})\varphi(q^{2})\varphi(q^{6})
+16q^{3}\varphi(q^{2})\psi(q^{12})\psi(q^{12})\psi(q^{4})
\\&\q+16q^{2}\psi(q^{4})\psi(q^{12})\psi(q^{4})\varphi(q^{6})
+16q^{2}\psi(q)\psi(q^{3})\psi(q^{12})\psi(q^{4}).
\endalign$$
Now applying (1.7) we get
$$\aligned&\sum_{n=0}^{\infty}N(1,3,4,12;8n+4)q^{n}
\\&=8\psi(q)\psi(q^{3})\varphi(q^{2})\varphi(q^{6})
+32q^{2}\psi(q)\psi(q^{3})\psi(q^{12})\psi(q^{4}).
\endaligned\tag 3.1$$
On the other hand, using (1.9) we see that
$$\align&\sum_{n=0}^{\infty}N(1,3,4,12;n)q^{n}
\\&=\big(\varphi(q^{4})+2q\psi(q^8)\big)\big(\varphi(q^{12})
+2q^3\psi(q^{24})\big)\varphi(q^{4})\varphi(q^{12})
\endalign$$
and so
$$\align&\sum_{n=0}^{\infty}N(1,3,4,12;2n+1)q^{2n+1}
\\&=2q\psi(q^8)\varphi(q^{12})^2\varphi(q^{4})
+2q^3\varphi(q^{4})\psi(q^{24})\varphi(q^{4})\varphi(q^{12}).
\endalign$$
Replacing $q$ with $q^{1/2}$ we then obtain
$$\align&\sum_{n=0}^{\infty}N(1,3,4,12;2n+1)q^{n}
\\&=2\psi(q^4)\varphi(q^{6})^2\varphi(q^{2})
+2q\varphi(q^{2})\psi(q^{12})\varphi(q^{2})\varphi(q^{6}).
\endalign$$
Now applying (1.7) we get
$$\sum_{n=0}^{\infty}N(1,3,4,12;2n+1)q^{n}
=2\psi(q)\psi(q^3)\varphi(q^{6})\varphi(q^{2}).\tag 3.2$$ From (3.1)
and (3.2) we deduce that
$$\align&\sum_{n=0}^{\infty}N(1,3,4,12;8n+4)q^{n}
-4\sum_{n=0}^{\infty}N(1,3,4,12;2n+1)q^{n}
\\&=32q^{2}\psi(q)\psi(q^{3})\psi(q^{12})\psi(q^{4})
=32q^{2}\sum_{n=0}^{\infty}t'(1,3,4,12;n)q^n
\\&=2q^{2}\sum_{n=0}^{\infty}t(1,3,4,12;n)q^n
.\endalign$$ Comparing the coefficients of $q^{n+2}$ on both sides
we obtain
$$t(1,3,4,12;n)=\f 12N(1,3,4,12;8n+20)-2N(1,3,4,12;2n+5).\tag 3.3$$
Now applying Lemma 3.4 we deduce the result.

\section*{4. Formulas for $t(1,1,3,4;n)$, $t(1,1,5,5;n)$,
 $t(1,5,5,5;n)$,
$t(1,3,3,12;n)$, $t(1,1,1,12;n)$, $t(1,1,3,12;n)$ and
$t(1,3,3,4;n)$}
\par For $a,b,c,d,n\in\Bbb N$ let
$$N_0(a,b,c,d;n)=\big|\big\{(x,y,z,w)\in\Bbb Z^4\bigm|
n=ax^2+by^2+cz^2+dw^2,\ 2\nmid xyzw\big\}\big|.$$ From [WS, (4.1)]
we know that
$$t(a,b,c,d;n)=N_0(a,b,c,d;8n+a+b+c+d).\tag 4.1$$
\par For $n\in \Bbb N$ following [AALW4] we define
$$\aligned &A(n)=\sum_{d\mid n}d\Ls{12}{n/d},
\q B(n)=\sum_{d\mid n}d\Ls{-3}{d}\Ls{-4}{n/d},
\\&C(n)=\sum_{d\mid n}d\Ls{-3}{n/d}\Ls{-4}{d},
\q D(n)=\sum_{d\mid n}d\Ls{12}{d},
\\&E(n)=\sum\Sb (i,j)\in\Bbb N\times \Bbb N\\i,j\ odd\\4n=i^2+3j^2\endSb
(-1)^\f{i-1}2 i\qtq {and} F(n)=\sum\Sb (i,j)\in\Bbb N\times\Bbb
N\\i,j\ odd\\4n=i^2+3j^2\endSb (-1)^\f{j-1}2 j.\endaligned$$ Suppose
that $n=2^{\a}3^{\beta}n_1$, where $\a$ and $\beta$ are non-negative
integers, $n_1\in\Bbb N$ and
 $\t{gcd}(n_1,6)=1$. From [AALW4, Theorem 3.1] we know that
 $$\aligned
&A(n)=2^{\a}3^{\beta}A(n_1), \q
B(n)=(-1)^{\a+\beta}2^{\a}\Ls{-3}{n_1}A(n_1),
\\&C(n)=(-1)^{\a+\beta+\f{n_1-1}2}3^{\beta}A(n_1)
\qtq{and} D(n)=\Ls 3{n_1}A(n_1).\endaligned\tag 4.2$$

 \par\pro{Lemma 4.1 ([AALW3, Theorem 7.2])} Let $n\in\Bbb N$ with $n\e1\mod2$. Then
$$N(1,1,3,4;n)=3A(n)-B(n)+\f 32C(n)-\f 12D(n)+E(n).$$
\endpro
\pro{Lemma 4.2 ([AALW3, Theorem 7.2])} Let $n\in\Bbb N$ with
$n\e1\mod2$. Then
$$N(1,1,4,12;n)=\f 32A(n)-\f 12B(n)+\f 32C(n)-\f 12D(n)+
\f 12E(n)+\f 32F(n).$$
\endpro
\pro{Theorem 4.1} Suppose $n\in\Bbb N$ and $8n+9=3^{\beta}n_1$ with
$3\nmid n_1$. Then
$$ t(1,1,3,4;n)=\f 12\Big(3^{\beta+1}\Ls 3{n_1}-1\Big)
\sum_{d\mid n_1}d\Ls 3d -\sum\Sb a,b\in\Bbb N,\ 2\nmid
a\\4(8n+9)=a^2+3b^2\endSb (-1)^{\f{a-1}2}a.$$\endpro

Proof. Since
$$\align &N(1,1,3,4;8n+9)\\&=\big|\{(x,y,z,w)\in \Bbb Z^4\ |\ 8n+9=x^2+y^2+3z^2+4w^2
\}\big|\\&=\big|\{(x,y,z,w)\in \Bbb Z^4\ |\
8n+9=x^2+y^2+3(2z)^2+4w^2 \}\big|\\&\q+\big|\{(x,y,z,w)\in \Bbb Z^4\
|\ 8n+9=x^2+y^2+3z^2+4w^2,2\nmid z \}\big|
\\&=\big|\{(x,y,z,w)\in \Bbb Z^4\ |\ 8n+9=x^2+y^2+12z^2+4w^2
\}\big|\\&\q+\big|\{(x,y,z,w)\in \Bbb Z^4\ |\
8n+9=x^2+y^2+3z^2+4w^2,2\nmid {xyzw} \}\big|
\\&=N(1,1,4,12;8n+9)+N_0(1,1,3,4;8n+9)
\\&=N(1,1,4,12;8n+9)+t(1,1,3,4;n),\endalign$$
we have $t(1,1,3,4;n)=N(1,1,3,4;8n+9)-N(1,1,4,12;8n+9).$ Now
applying Lemmas 4.1, 4.2 and (4.2) we deduce the result.

\par{\bf Remark 4.1} Theorem 4.1 was conjectured by the authors in
[WS].
  \pro{Lemma 4.3 ([AAW, Theorem 7.1])} Let $n\in\Bbb N$ and
$n=2^{\a}5^{\beta}n_1$ with $n_1\in\Bbb N$ and  $\t{gcd}(n_1,10)=1$.
Then
$$\aligned N(1,1,5,5;n)=\cases 2(5^{\beta+1}-3)\sigma(n_1)
&\t{if $2\mid n$,}
\\\f 23(5^{\beta+1}-3)\sigma(n_1)+\f 83c(n)&\t{if $2\nmid n$.}
\endcases\endaligned$$
where $c(n)$ is given by
$$q\prod_{n=1}^{\infty}(1-q^{2n})^2(1-q^{10n})^2=\sum_{n=1}^{\infty}
 c(n)q^n.$$
\endpro
\pro{Theorem 4.2} Suppose $n\in\Bbb N$ and $2n+3=5^{\beta}n_1$ with
$5\nmid n_1$. Then
$$t(1,1,5,5;n)=\f 43(5^{\beta+1}-3)\sigma(n_1)-\f83c(2n+3).$$\endpro
Proof. Note that
$$\align &N(1,1,5,5;8n+12)\\&=\big|\{(x,y,z,w)\in \Bbb Z^4\ |\ 8n+12
=x^2+y^2+5z^2+5w^2 \}\big|\\&=\big|\{(x,y,z,w)\in \Bbb Z^4\ |\
8n+12=(2x)^2+(2y)^2+5(2z)^2+5(2w)^2 \}\big|\\&\q+\big|\{(x,y,z,w)\in
\Bbb Z^4\ |\ 8n+12=x^2+y^2+5z^2+5w^2,2\nmid {xyzw} \}\big|
\\&=N(1,1,5,5;2n+3)+N_0(1,1,5,5;8n+12)
\\&=N(1,1,5,5;2n+3)+t(1,1,5,5;n),\endalign$$
applying Lemma 4.3 we deduce the result.

 \pro{Lemma 4.4 ([AAW, Theorem 6.1])} Let $n\in\Bbb N$. Then
$$ N(1,5,5,5;n)=\sum_{d\mid n}(-1)^{n+d}\Big(\Ls 5d+\Ls 5{n/d}\Big)d.$$
\endpro

\pro{Theorem 4.3} Let $n\in\Bbb N$. Then
$$\align t(1,5,5,5;n)&=\sum_{d\mid {8n+16}}(-1)^{d}\Big(\Ls 5d +\Ls 5{{(8n+16)}/d}
\Big)d\\&\q-\sum_{d\mid {2n+4}}(-1)^{d}\Big(\Ls 5d+\Ls
5{{(2n+4)}/d}\Big) d\endalign$$\endpro
 Proof. Observe that
$$\align &N(1,5,5,5;8n+16)\\&=\big|\{(x,y,z,w)\in \Bbb Z^4\ |
\ 8n+16=x^2+5y^2+5z^2+5w^2 \}\big|\\&=\big|\{(x,y,z,w)\in \Bbb Z^4\
|\ 8n+16=(2x)^2+5(2y)^2+5(2z)^2+5(2w)^2
\}\big|\\&\q+\big|\{(x,y,z,w)\in \Bbb Z^4\ |\
8n+16=x^2+5y^2+5z^2+5w^2,2\nmid {xyzw} \}\big|
\\&=N(1,5,5,5;2n+4)+N_0(1,5,5,5;8n+16)
\\&=N(1,5,5,5;2n+4)+t(1,5,5,5;n),\endalign$$
applying Lemma 4.4 we deduce the result.
 \pro{Lemma 4.5 ([AALW3, Theorem
7.2])} Let $n\in\Bbb N$ with $2\nmid n$. Then
$$ N(1,3,3,12;n)=A(n)+B(n)-\f12C(n)-\f12D(n)+F(n)$$
and
$$ N(3,3,4,12;n)=\f12(A(n)+B(n)-C(n)-D(n)-E(n)+F(n)).$$
\endpro
\pro{Theorem 4.4} Let $n\in\Bbb N$ and $8n+19=3^{\beta}n_1$ with
$n_1\in\Bbb N$ and $3\nmid n_1$. Then
$$\aligned &t(1,3,3,12;n)\\&=\f 12\Big(3^{\beta}\Ls 3{n_1}-1\Big)
\sum_{d\mid n_1}d\Ls 3d+\f12\sum\Sb a,b\in\Bbb N,\ a\e b \e 1\mod
2\\4(8n+19)=a^2+3b^2\endSb
((-1)^{\f{a-1}2}a+(-1)^{\f{b-1}2}b).\endaligned$$
\endpro
Proof.
Observe that
$$\align &N(1,3,3,12;8n+19)\\&=\big|\{(x,y,z,w)\in \Bbb Z^4\ |\ 8n+19=x^2+3y^2+3z^2+12w^2
\}\big|\\&=\big|\{(x,y,z,w)\in \Bbb Z^4\ |\
8n+19=(2x)^2+3y^2+3z^2+12w^2 \}\big|\\&\q+\big|\{(x,y,z,w)\in \Bbb
Z^4\ |\ 8n+19=x^2+3y^2+3z^2+12w^2,2\nmid x \}\big|
\\&=\big|\{(x,y,z,w)\in \Bbb Z^4\ |\ 8n+19=4x^2+3y^2+3z^2+12w^2
\}\big|\\&\q+\big|\{(x,y,z,w)\in \Bbb Z^4\ |\
8n+19=x^2+3y^2+3z^2+12w^2,2\nmid {xyzw} \}\big|
\\&=N(4,3,3,12;8n+19)+N_0(1,3,3,12;8n+19)
\\&=N(3,3,4,12;8n+19)+t(1,3,3,12;n),\endalign$$
applying Lemma 4.5 and (4.2) we obtain the result.
 \pro{Lemma 4.6 ([AALW3,
Theorem 7.2])} Let $n\in\Bbb N$ and $2\nmid n$. Then
$$ N(1,1,1,12;n)=3A(n)-B(n)+\f32C(n)-\f12D(n)+3F(n).$$
\endpro
\pro{Theorem 4.5} Let $n\in\Bbb N$ and $8n+15=3^{\beta}n_1$ with
$n_1\in\Bbb N$ and $3\nmid n_1$. Then
$$ t(1,1,1,12;n)=\f 12\Big(3^{\beta+1}\Ls 3{n_1}+1\Big)
\sum_{d\mid n_1}d\Ls 3d+3\sum\Sb a,b\in\Bbb N,\ 2\nmid
a\\4(8n+15)=a^2+3b^2\endSb (-1)^{\f{b-1}2}b.$$\endpro
 Proof. Since $8n+15=x^2+y^2+z^2+12w^2$ for $x,y,z,w\in\Bbb Z$
 implies that $2\nmid xyzw$, we see that
 $$t(1,1,1,12;n)=N_0(1,1,1,12;8n+15)=N(1,1,1,12;8n+15).$$
 Now the result follows from Lemma 4.6 and (4.2).

 \pro{Lemma 4.7 ([AALW1, Theorems 1.10 and 1.13])} Let
$n\in\Bbb N$ and $n=3^{\beta}n_1$ with $3\nmid n_1$. For $n\e1\mod4$
we have
$$N(1,1,3,12;n)=3\sigma(n_1)+a(n)\qtq{and}
N(1,1,12,12;n)=2\sigma(n_1)+2a(n).$$
\endpro
\pro{Theorem 4.6} Suppose $n\in\Bbb N$ and $8n+17=3^{\beta}n_1$ with
$3\nmid n_1$. Then
$$ t(1,1,3,12;n)=\sigma(n_1)-a(8n+17).$$\endpro
Proof. Since
$$\align &N(1,1,3,12;8n+17)\\&=\big|\{(x,y,z,w)\in \Bbb Z^4\ |\ 8n+17=x^2+y^2+3z^2+12w^2
\}\big|\\&=\big|\{(x,y,z,w)\in \Bbb Z^4\ |\
8n+17=x^2+y^2+3(2z)^2+12w^2 \}\big|\\&\q+\big|\{(x,y,z,w)\in \Bbb
Z^4\ |\ 8n+17=x^2+y^2+3z^2+12w^2,2\nmid z \}\big|
\\&=\big|\{(x,y,z,w)\in \Bbb Z^4\ |\ 8n+17=x^2+y^2+12z^2+12w^2
\}\big|\\&\q+\big|\{(x,y,z,w)\in \Bbb Z^4\ |\
8n+17=x^2+y^2+3z^2+12w^2,2\nmid {xyzw} \}\big|
\\&=N(1,1,12,12;8n+17)+N_0(1,1,3,12;8n+17)
\\&=N(1,1,12,12;8n+17)+t(1,1,3,12;n),\endalign$$
 applying
Lemma 4.7 we deduce the result.
 \pro{Lemma 4.8 ([AALW1,
Theorems 1.16 and 1.23])} Let $n\in\Bbb N$ and $n=3^{\beta}n_1$ with
$3\nmid n_1$. For
 $n\e3\mod4$ we have
$$N(1,3,3,4;n)=3\sigma(n_1)-a(n)\qtq{and}
N(3,3,4,4;n)=2\sigma(n_1)-2a(n).$$
\endpro
\pro{Theorem 4.7} Suppose $n\in\Bbb N$ and $8n+11=3^{\beta}n_1$ with
$3\nmid n_1$. Then
$$ t(1,3,3,4;n)=\sigma(n_1)+a(8n+11).$$\endpro
Proof. Since
$$\align &N(1,3,3,4;8n+11)\\&=\big|\{(x,y,z,w)\in \Bbb Z^4\ |\ 8n+11=x^2+3y^2+3z^2+4w^2
\}\big|\\&=\big|\{(x,y,z,w)\in \Bbb Z^4\ |\
8n+11=(2x)^2+3y^2+3z^2+4w^2 \}\big|\\&\q+\big|\{(x,y,z,w)\in \Bbb
Z^4\ |\ 8n+11=x^2+3y^2+3z^2+4w^2,2\nmid x \}\big|
\\&=\big|\{(x,y,z,w)\in \Bbb Z^4\ |\ 8n+11=4x^2+3y^2+3z^2+4w^2
\}\big|\\&\q+\big|\{(x,y,z,w)\in \Bbb Z^4\ |\
8n+11=x^2+3y^2+3z^2+4w^2,2\nmid {xyzw} \}\big|
\\&=N(3,3,4,4;8n+11)+N_0(1,3,3,4;8n+11)
\\&=N(3,3,4,4;8n+11)+t(1,3,3,4;n),\endalign$$
applying Lemma 4.8 we derive the result.

\end{document}